\documentclass[11pt]{article}

\usepackage[a4paper,margin=1in]{geometry}
\usepackage{amsmath,amssymb,amsthm}
\usepackage{graphicx}
\usepackage[colorlinks=true,citecolor=black,linkcolor=black,urlcolor=blue]{hyperref}

\theoremstyle{plain}
\newtheorem{theorem}{Theorem}
\newtheorem{lemma}[theorem]{Lemma}
\newtheorem{corollary}[theorem]{Corollary}
\newtheorem{proposition}[theorem]{Proposition}

\theoremstyle{definition}

\newtheorem{conjecture}[theorem]{Conjecture}

\theoremstyle{remark}

\title{$k$-ary Legendre pairs}

\author{Ilias S. Kotsireas\\
\small CARGO Lab, Wilfrid Laurier University\\
\small Waterloo, Ontario, Canada\\
\small \texttt{ikotsire@gmail.com}
\and
Huaning Liu\\
\small School of Mathematics, Northwest University\\
\small Xi'an, China\\
\small \texttt{hnliu@nwu.edu.cn}
\and
Arne Winterhof\\
\small RICAM, Austrian Academy of Sciences\\
\small Linz, Austria\\
\small \texttt{arne.winterhof@oeaw.ac.at}}

\date{}

\begin{document}

	\maketitle

	\begin{abstract}
		We study $k$-ary Legendre pairs and derive arithmetic
		restrictions on their existence and balance. For an odd prime
		$p$, we determine the component sums of a $p$-ary Legendre pair
		and show that, after normalization, each component has fixed
		symbol composition. As a consequence, every $p$-ary Legendre
		pair of length $\ell$ yields a Butson Hadamard matrix
		$H(2p,2\ell+2)$.
		
		For ternary Legendre pairs we combine fixed-content necklaces,
		integral-PAF symmetries, exact half-frequency restrictions,
		compression-decompression, and PSD tests. These reductions lead
		to ternary Legendre pairs for every admissible length
		$\ell\le32$. We conjecture that ternary Legendre pairs exist for every admissible length.
	\end{abstract}

\medskip
\noindent\textbf{Keywords.} Legendre pairs; periodic autocorrelation; Butson Hadamard matrices; ternary sequences; power spectral density; sequence design.

\smallskip
\noindent\textbf{2020 Mathematics Subject Classification.} Primary: 94A55; Secondary: 05B20, 11T71.


	\section{Introduction}
	
	\quad\ \ Let $A=[a_0, a_1, \ldots, a_{\ell-1}]\in\mathbb{C}^{\ell}$ be a finite complex
	sequence of length $\ell\ge 2$ periodically continued with period $\ell$. 
	The {\em periodic autocorrelation function} $PAF(A, s)$ of $A$ at lag $s$ is
	$$
	PAF(A,s)=\sum_{j=0}^{\ell-1}a_j\overline{a_{j+s}}, \qquad s=0,1,\ldots, \ell-1,
	$$
	where $\overline{b}$ is the complex conjugate of $b\in \mathbb{C}$. Two sequences
	$A=[a_0, a_1, \ldots, a_{\ell-1}]$ and $B=[b_0, b_1, \ldots, b_{\ell-1}]$ of length $\ell$ form
	a {\em Legendre pair} $(A, B)$ if
	\begin{eqnarray}\label{Equation:KW-1}
		PAF(A,s)+PAF(B,s)=-2, \qquad s=1,2,\ldots,\left\lfloor\frac{\ell}{2}\right\rfloor.
	\end{eqnarray}
	Since $PAF(A,\ell-s)=\overline{PAF(A,s)}$ for $s=1,2,\ldots,\left\lfloor\frac{\ell}{2}\right\rfloor$ 
	we may extend the range of $s$ in (\ref{Equation:KW-1}) to $s=1,2,\ldots, \ell-1$.
	
	A Legendre pair $(A, B)$ is {\em binary} if $A, B\in\{-1,+1\}^{\ell}$.
	A {\em Hadamard matrix} is an $n\times n$ matrix $H$ with elements $\pm 1$ such that $HH^{T}=nI_n$,
	where $I_n$ is the $n\times n$ identity matrix.	
	Binary Legendre pairs of length~$\ell$ are pertinent to the construction of
	Hadamard matrices of order $2\ell+2$, see for example \cite[Theorem 2]{ArasuBH2020}, \cite[Theorem 3]{FletcherGS2001} or
	\cite[Theorem 1]{Kotsireas2013}. Binary Legendre pairs of length $\ell$ are known for several infinite sets of $\ell$,
	if $\ell>2$ is a prime (Legendre sequences),
	if $\ell$ is odd and $2\ell+1$ is a prime power (Szekeres), 
	if $\ell=2^{r}-1$ with $r\geq 2$, ($m$-sequences),
	and a set of $\ell$ which is infinite if the twin-prime conjecture is true,
	if $\ell=p(p+2)$ where $p$ and $p+2$ are twin
	primes, (Jacobi sequence or twin-prime sequence).
	Several ad hoc constructions are also known. The smallest integer for which the question of existence of binary
	Legendre pairs is still undecided is $\ell=115$.
	See for example \cite{KotsireasW2024} and references therein for the results mentioned here.
	
	A Legendre pair $(A, B)$ is {\em quaternary} if $A, B\in\{1, i, -1, -i\}^{\ell}$,
	where $i=\sqrt{-1}$. The first and the third author \cite{KotsireasW2024} introduced
	quaternary Legendre pairs of length $\ell$ and showed that
	they are pertinent to the construction of {\em quaternary Hadamard matrices} of
	order $2\ell+2$ and thus of binary Hadamard matrices of order $4\ell+4$. 
	Since a binary Legendre pair can be considered quaternary we may restrict ourselves to quaternary pairs of even length $\ell$. 
	Jedwab and Pender \cite{JP2025} showed that 
	there is a quaternary Legendre pair of length $\ell$
	if $2\ell+1$ is a prime power and
	if $\ell=2p$ for an odd prime $p$ for which $2p-1$ is a prime power. 
	The first author, Koutschan and the third author \cite{KotsireasKW2025,KotsireasW2024}
	found the remaining quaternary Legendre pairs of lengths $\ell\le 40$, Lebedev \cite{Lebedev2026} completed the list finding $7$ more values and Verma \cite{v26} settled $\ell=64$. Now $\ell=70,76,88,92,94$ and $100$ are the smallest undecided cases. 
	
	
	A Legendre pair $(A, B)$ is {\em $k$-ary} if $A, B\in\{1, \xi_k, \xi_k^2, \ldots, \xi_k^{k-1}\}^{\ell}$,
	where $\xi_k=\hbox{e}^{2\pi \sqrt{-1}/k}$ is a primitive complex $k$-th root of unity. For a complex matrix $M$, let $M^*=\overline{M}^{\,T}$ denote its conjugate transpose. A {\em Butson Hadamard matrix} $H(k,n)$ is an $n\times n$-matrix with entries in $\{1, \xi_k, \xi_k^2, \ldots, \xi_k^{k-1}\}$ which satisfies
	$H(k,n)H(k,n)^{*}=nI_n$. From~\cite{Butson1962,Winterhof2000} we know that if $p$ is prime, $f\ge 1$ and $n>1$, then $H(p^f,n)$ can exist only for $p\mid n$.
    For further non-existence results see for example \cite{Brock1988,Winterhof2000}.
	
	Our main contributions of this paper are the following results.
	
	First, in Section \ref{non-ex} we prove that for prime $k$ a $k$-ary Legendre pair of length $\ell$ can only exist if $\ell\equiv -1\bmod k.$ 	
	
	Second, in Section \ref{Section:Butson Hadamard matrices} we study the connection between Legendre pairs and
	Butson Hadamard matrices 
	for general $k$-ary Legendre pairs. More precisely, under
	explicit conditions on the sequence sums of a Legendre pair of length $\ell$, we construct Butson
	Hadamard matrices $H(2k,2\ell+2)$ when $k$ is odd and
	$H(k,2\ell+2)$ when $k$ is even. 
	For odd primes $k=p$, 
    these conditions are automatically satisfied, yielding a Butson Hadamard matrix $H(2p,2\ell+2)$ from every $p$-ary Legendre pair of length $\ell$.
    Note that a $H(p,2\ell+2)$ may not exist even if a $p$-ary Legendre pair of length $\ell\equiv -1\bmod p$ exists. 
	
	Third, Section \ref{ternary} deals with the case $k=3$, that is, ternary Legendre pairs. 
	A pair of (generalized) $m$-sequences covers all lengths of the form $\ell=3^r-1$ with $r\geq1$.
	We combine fixed-content necklaces, integral-PAF symmetries, and exact half-frequency restrictions with compression-decompression and PSD tests. In this way we obtain ternary Legendre pairs for every admissible length $\ell\leq 32$. The smallest open cases are then $\ell=35$ and $\ell=38$.
    Note that for $\ell\equiv -1\bmod 3$ and $\ell\le 32$ there is no $BH(3,2\ell+2)$ by \cite[Theorem 5]{Winterhof2000} up to the possible exceptions $\ell=5,17$ and $23$.
    However, as we proved in Section \ref{Section:Butson Hadamard matrices}, for all these $\ell$ there is a $BH(6,2\ell+2)$.
	
	\section{Non-existence, balance and a construction}	\label{non-ex}
	
	\subsection{A non-existence result}
	
	In this subsection we derive a necessary condition on the length $\ell$ for the existence of a $k$-ary Legendre pair when $k$ is prime.
	
	\begin{theorem}\label{Theorem:length}
		For prime $k$, a $k$-ary Legendre pair of length $\ell\ge2$ can exist only if $\ell\equiv -1\pmod k$.
	\end{theorem}
	\begin{proof}
		For $k=2$, the assertion is well-known. 
		It follows from 
		$PAF(A,s)\equiv PAF(B,s)\equiv \ell\bmod 4,$
		see for example \cite[Proposition 6.4.16]{niwi15}, and \eqref{Equation:KW-1}.
		
		Now we assume that $k$ is an odd prime. 
		Let $(A, B)$ be a $k$-ary Legendre pair of length $\ell$.  For any $s=1,\ldots,\ell-1$, we have
		\begin{eqnarray*}
			0=PAF(A, s)+PAF(B, s)+2=\sum_{j=0}^{\ell-1}a_j\overline{a_{j+s}}
			+\sum_{j=0}^{\ell-1}b_j\overline{b_{j+s}}+2
			=\sum_{j=0}^{k-1}c_j \xi_k^j,
		\end{eqnarray*}
		where $c_j$ counts the occurrences of $\xi_k^j$  
		in the two autocorrelation sums and two additional $1$s.
		Thus the $c_j$ are nonnegative integers satisfying  $\sum\limits_{j=0}^{k-1}c_j=2\ell+2$. 
		Since
		$\sum\limits_{j=0}^{k-1}\xi_k^j=0$, we get
		$0=\sum\limits_{j=0}^{k-2}(c_j-c_{k-1})\xi_k^j.
		$ 
		As $k$ is prime, the cyclotomic polynomial $\Phi_k(X)=1+X+\cdots+X^{k-1}$ has degree $k-1$ and is irreducible over $\mathbb{Q}$. Hence $1,\xi_k,\ldots,\xi_k^{k-2}$ are linearly independent over $\mathbb{Q}$. Therefore $c_0=\ldots=c_{k-1}$ and thus 
		$2\ell+2=kc_0\equiv 0\bmod k,$ that is,
		\begin{eqnarray}\label{Equation: Non-existence result}
			\ell\equiv -1\bmod k,
		\end{eqnarray}
		and the result follows. 	
	\end{proof}
	
	
	\subsection{Balance}
	
	\quad\ \ In this subsection we study the sequence sums of the two
	components of a $k$-ary Legendre pair. We first establish a
	condition on the balance for general $k$, and then determine the
	sequence sums and symbol compositions completely when the
	alphabet size is an odd prime. 

	Write $A=[a_0, a_1, \ldots, a_{\ell-1}]$, $B=[b_0, b_1, \ldots, b_{\ell-1}]$
	and let $(A,B)$ be a $k$-ary Legendre pair of length~$\ell$. Put
	\begin{eqnarray}\label{Definition:alpha+beta}
		\alpha=\sum_{j=0}^{\ell-1}a_j \qquad \hbox{and} \qquad \beta=\sum_{j=0}^{\ell-1}b_j.
	\end{eqnarray}
	
	\begin{lemma}\label{Lemma:k-ary-0}
		For any $k$-ary Legendre pair and $\alpha,\beta$ defined by $(\ref{Definition:alpha+beta})$
		we have
		\begin{eqnarray}\label{Equation:alpha2+beta2=2}
			|\alpha|^2+|\beta|^2=2.
		\end{eqnarray}
	\end{lemma}
	
	\begin{proof}  
		For $j=0,1,\ldots,\ell-1$ we have $\sum\limits_{s=0}^{\ell-1}\overline{a_{j+s}} =\sum\limits_{s=0}^{\ell-1}\overline{a_s}$. 		
		Hence,
		$$
		\sum_{s=0}^{\ell-1}PAF(A, s)=\sum_{s=0}^{\ell-1}\sum_{j=0}^{\ell-1}a_{j}\overline{a_{j+s}}
		=\sum_{j=0}^{\ell-1}a_{j}\sum_{s=0}^{\ell-1}\overline{a_{s}}=\left|\sum_{j=0}^{\ell-1}a_{j}\right|^2
		=|\alpha|^2
		$$		
		and similarly $\sum\limits_{s=0}^{\ell-1}PAF(B, s)=|\beta|^2.
		$
		Therefore
		\begin{eqnarray*}
			|\alpha|^2+|\beta|^2=\sum_{s=0}^{\ell-1}\left(PAF(A, s)+PAF(B, s)\right)
			=2\ell+(\ell-1)\cdot(-2)=2
		\end{eqnarray*}
		by \eqref{Equation:KW-1}.
	\end{proof}
	
	Multiplying all entries of either component by the same $k$th root of unity
    preserves its PAFs. Independent cyclic shifts of the two components also preserve
    their PAFs. Simultaneously conjugating both components preserves
    condition~\eqref{Equation:KW-1}, as does interchanging the two components.
    Conjugating only one component need not preserve the Legendre pair condition.
	Hence, for binary Legendre pairs we may assume 
	$\alpha=\beta=1$.
	For quaternary Legendre pairs we may assume
	$\alpha=0,$ $\beta=1+i$ for even $\ell$
	and
	$\alpha=\beta=1$ for odd $\ell$,
	see \cite{KotsireasW2024}.
	
	\begin{lemma}\label{lem:balance-rigidity}
		Let $(A,B)$ be a $k$-ary Legendre pair. Then
		$$
		(|\alpha|^2,|\beta|^2)\in\{(0,2),(1,1),(2,0)\}.
		$$
	\end{lemma}
	
	\begin{proof}
		Put $x=|\alpha|^2$. Since $\alpha\in\mathbb Z[\xi_k]$ 
		we have $x\in \mathbb Z[\xi_k]\cap \mathbb R=\mathbb Z[\xi_k+\xi_k^{-1}]$. 

        An automorphisms $\sigma_t$ 
		of $\mathbb Q(\xi_k)$ is uniquely defined by $\sigma_t(\xi_k)=\xi_k^t$ with $\gcd(t,k)=1$, $t=1,\ldots,k-1$. 
		If $\sigma_t$ is applied entrywise to $A$ and $B$, then
		for every lag $s$ we have
		$$
		PAF(\sigma_t(A),s)=\sum_{j=0}^{\ell-1}\sigma_t(a_j)\overline{\sigma_t(a_{j+s})}
		=\sum_{j=0}^{\ell-1}\sigma_t(a_j) \sigma_t(\overline{a_{j+s}})=\sigma_t(PAF(A,s))
		$$
		as well as $PAF(\sigma_t(B),s)=\sigma_t(PAF(B,s))$. 
        Hence, $(\sigma_t(A),\sigma_t(B))$ is again a
		$k$-ary Legendre pair.

        By Lemma~\ref{Lemma:k-ary-0} we get $|\sigma_t(\alpha)|^2+|\sigma_t(\beta)|^2=2$.
		In particular, we have
		$$
		0\le \sigma_t(x)=\sigma_t(\alpha\overline{\alpha})=|\sigma_t(\alpha)|^2\le 2.
		$$

        Put $y=x-1$.
        Note that $y\in \mathbb Z[\xi_k+\xi_k^{-1}]$ is an algebraic integer and thus its norm, the product of its conjugates $\sigma_t(y)=\sigma_t(x)-1$, $\gcd(t,k)=1$, $t=0,\ldots,\lfloor k/2\rfloor$, is a rational integer in  
		$[-1,1]$, that is, $y\in \{-1,0,1\}$ and thus $x=|\alpha|^2\in \{0,1,2\}$. 
        Finally,
		Lemma~\ref{Lemma:k-ary-0} gives
		$|\alpha|^2+|\beta|^2=2$, and the result follows. 
	\end{proof}

     Let $c_r$ be the number of $n=0,1,\ldots,\ell-1$ with $a_n=\xi_k^r$. Then we call
    $(c_0,c_1,\ldots,c_{k-1})$ the {\em symbol composition} of $A=[a_0,a_1,\ldots,a_{\ell-1}]$.
	
	\begin{theorem}	\label{thm:prime-balance}
		Let $p$ be an odd prime and let $(A,B)$ be a $p$-ary
		Legendre pair of length $\ell\ge2$. Then
		$$
		\alpha,\beta \in \{-1,-\xi_p,\ldots,-\xi_p^{p-1}\}.
		$$
		Consequently, after multiplying $A$ and $B$ independently by
		suitable $p$th roots of unity, we may assume
		$$
		\alpha=\beta=-1.
		$$
		Moreover, writing $\ell=ph-1$, after this normalization
		each component has symbol composition
		$(c_0,c_1,\ldots,c_{p-1})=(h-1,h,\ldots,h)$.         
\end{theorem}
	
	\begin{proof}
		We first exclude the possibility that one of the sequence sums
		vanishes. 
		Obviously, we have
        $$
		\sum_{r=0}^{p-1}c_r=\ell\quad \mbox{and} 	\quad 	\alpha=\sum_{r=0}^{p-1}c_r\xi_p^r.
		$$
		Suppose that $\alpha=0$ and define
		$P(X)=\sum_{r=0}^{p-1}c_rX^r\in\mathbb Z[X]$.
		Then $P(\xi_p)=0$. Since $p$ is prime, $\Phi_p(X)=1+X+\cdots+X^{p-1}$
		is the minimal polynomial of $\xi_p$ over $\mathbb Q$, of degree
		$p-1$. Hence $\Phi_p$ divides $P$ in $\mathbb Q[X]$. Because
		$\deg P\le p-1$, we must have $P(X)=c\,\Phi_p(X)$ for some $c\in\mathbb Q$. 
		Comparing coefficients gives
		$c_0=c_1=\cdots=c_{p-1}=c$. Since the $c_r$ are integers,
		$c\in\mathbb Z$, and therefore $p\mid\ell$, contradicting
		Theorem~\ref{Theorem:length}, which gives
		$\ell\equiv-1\pmod p$. Thus $\alpha\ne0$. The same argument
		shows that $\beta\ne0$.

		Lemma~\ref{lem:balance-rigidity} now implies
		$|\alpha|^2=|\beta|^2=1$.
		
		We next determine the multiplicities $c_r$. For
		$0\le t\le p-1$, set
		$S_t=\sum\limits_{r=0}^{p-1}c_r\xi_p^{tr}$.
		For $t=0$ we have $S_0=\ell$. For $1\le t\le p-1$, the map
		$$
			\sigma_t:\xi_p\longmapsto\xi_p^t
		$$
		defines an automorphism of $\mathbb Q(\xi_p)$, and
		$S_t=\sigma_t(\alpha)$. 
        We have
		$$	
		|S_t|^2=\sigma_t(\alpha)\overline{\sigma_t(\alpha)}=\sigma_t(\alpha\overline{\alpha})=\sigma_t(1)=1.
		$$
		Thus
        $$\sum\limits_{t=0}^{p-1}|S_t|^2=\ell^2+p-1.$$

        Using the orthogonality
		relation
		$$
		\sum_{t=0}^{p-1}\xi_p^{t(r-s)}=\begin{cases}p,&r=s,\\
			0,&r\ne s,
			\end{cases}
		$$
		we obtain
		$$
		\sum_{t=0}^{p-1}|S_t|^2=\sum_{t=0}^{p-1} \sum_{r=0}^{p-1}\sum_{s=0}^{p-1} c_rc_s\xi_p^{t(r-s)}=p\sum_{r=0}^{p-1}c_r^2.
		$$
		
		Consequently, $p\sum\limits_{r=0}^{p-1}c_r^2=\ell^2+p-1$.
		By Theorem~\ref{Theorem:length}, we can write
		$\ell=ph-1$ for some positive integer $h$. Then 
        $$\sum\limits_{r=0}^{p-1}c_r=ph-1\quad \mbox{and} 
		\quad
		\sum_{r=0}^{p-1}c_r^2
        =ph^2-2h+1=(h-1)^2+(p-1)h^2.
		$$

        Obviously, $c_0=h-1$ and $c_r=h$ for $r=1,2,\ldots,p-1$ is a solution of this system of equations.
        Among
		all $p$-tuples of nonnegative integers with fixed sum $ph-1$, the
		sum of squares is minimized when the entries differ by at most
		one. Indeed, if $c_r\ge c_s+2$, replacing $(c_r,c_s)$ by
		$(c_r-1,c_s+1)$ preserves the total sum and changes the sum of
		squares by
		$$
		(c_r-1)^2+(c_s+1)^2-c_r^2-c_s^2=-2(c_r-c_s)+2<0.
		$$
		Thus a solution tuple must have entries differing by at most
		one, that is, a cyclic shift of $(h-1,h,h,\ldots,h)$. 
		Therefore
		$$
		\alpha=(h-1)\xi_p^a +h\sum_{\substack{0\le r\le p-1\\r\ne a}}\xi_p^r=h\sum_{r=0}^{p-1}\xi_p^r-\xi_p^a
			=-\xi_p^a.
		$$
		The same argument applied to $B$ gives
		$\beta=-\xi_p^b$ for some $b$.

		Finally, multiply every entry of $A$ by $\xi_p^{-a}$ and every
		entry of $B$ by $\xi_p^{-b}$ 
        and 
        we may assume $\alpha=\beta=-1$. 
	\end{proof}
		
	\subsection{$m$-sequences: $k=p$ is prime and $\ell=p^{r}-1$} \label{Subsection:m-sequences}

    We recall the following standard construction for the convenience of the reader.
    
	Let $k=p$ be prime, let $r\ge1$ with $p^r>2$, and put $\ell=p^r-1$. Let $g$ be a primitive
	element of $\mathbb F_{p^r}$, and   let
	$$
	\psi(x)=\exp\left(\frac{2\pi\sqrt{-1}}{p}
	\operatorname{Tr}_{\mathbb F_{p^r}/\mathbb F_p}(x)\right)
	$$  
	be the canonical additive character of $\mathbb F_{p^r}$. Define
	$$
	A=B=[a_0,a_1,\ldots,a_{\ell-1}],\qquad
	a_j=b_j=\psi(g^j-1),\quad 0\le j\le \ell-1.
	$$
	In particular, we have $a_0=b_0=\psi(0)=1$. For
	$1\le s\le\ell-1=p^r-2$, we have
	\begin{eqnarray*}
		PAF(A,s)&=&PAF(B,s)=\sum_{j=0}^{\ell-1}
		\psi(g^j-1)\overline{\psi(g^{j+s}-1)}=\sum_{j=0}^{\ell-1}\psi\bigl((1-g^s)g^j\bigr)\\
		&=&\sum_{\gamma\in\mathbb F_{p^r}^{*}}\psi\bigl((1-g^s)\gamma\bigr)=\sum_{\eta\in\mathbb F_{p^r}^{*}}\psi(\eta)
		=-1.
	\end{eqnarray*}
	Indeed, $g^s\ne1$ for $1\le s\le p^r-2$, so the map
	$$
	\gamma\longmapsto(1-g^s)\gamma
	$$
	is a bijection of $\mathbb F_{p^r}^{*}$. The last equality follows
	from the orthogonality of nontrivial additive characters; see, for
	example, \cite[Chapter 5.1]{LidlN1997}. Hence
	$$
	PAF(A,s)+PAF(B,s)=-2,\qquad 1\le s\le\ell-1,
	$$
	and therefore $(A,B)$ is a $p$-ary Legendre pair of length
	$\ell=p^r-1$. Note that the definition of Legendre pairs includes the case $A=B$, that is, $PAF(A,s)=-1$ for $1\le s\le \ell-1$.
	
	\section{$k$-ary Legendre pairs and Butson Hadamard matrices}\label{Section:Butson Hadamard matrices}
	
	The circulant matrix $C(A)$ defined by $A=\left[a_0, a_1, \ldots, a_{\ell-1}\right]$ is 
	$$
	C(A)=\left[\begin{array}{ccccc}
		a_0          &      a_1    &  \ldots   &   a_{\ell-2}    &   a_{\ell-1}      \\
		a_{\ell-1}   &      a_0    &  \ldots   &   a_{\ell-3}    &   a_{\ell-2}      \\
		\vdots       &  \vdots     &           & \vdots          &   \vdots          \\
		a_2          &      a_3    &  \ldots   &   a_{0}         &   a_{1}     \\
		a_1          &      a_2    &  \ldots   &   a_{\ell-1}    &   a_{0}
	\end{array}\right].
	$$
	Clearly, we have
	\begin{equation} \label{Equation:Sec2-CACB}
		C(A)C(B)=C(B)C(A),
	\end{equation}
	see for example \cite[Theorem 3.2.4]{D79},
	and
	\begin{eqnarray} \label{Equation:Sec2-CACA+CBCB}
		\left[\begin{array}{ll}
			C(A) & C(B)
		\end{array}
		\right]
		\left[\begin{array}{l}
			C(\overline{A})^{T} \\
			C(\overline{B})^{T}
		\end{array}
		\right]
		&=&\left[PAF(A, u-v)+PAF(B, u-v)\right]_{1\leq u, v\leq \ell} \nonumber \\
		&=&\left[\begin{array}{llll}
			2\ell     &    -2     &  \ldots   &  -2   \\
			-2        &   2\ell   &  \ldots   &  -2   \\
			\vdots    & \vdots    &           & \vdots \\
			-2        &   -2      &  \ldots   &  2\ell
		\end{array}
		\right]_{\ell\times\ell}.
	\end{eqnarray}
	
	To border the circulant block matrix, we seek two additional coordinates making the first two rows orthogonal to the lower block rows. The resulting linear equations lead naturally to the following choices of $z, w$ (and $x, y$ in the even case):
	\begin{eqnarray}
		&&z=\frac{1}{2}\left(\alpha-\beta-\overline{\alpha}-\overline{\beta}\right), \qquad
		w=\frac{1}{2}\overline{\xi_{2k}}\left(\beta-\alpha-\overline{\alpha}-\overline{\beta}\right), \label{Definition:zw} \\
		&&x=\frac{1}{2}\left(\alpha-\beta-\overline{\alpha}-\overline{\beta}\right), \qquad
		y=\frac{1}{2}\left(\beta-\alpha-\overline{\alpha}-\overline{\beta}\right). \label{Definition:xy}
	\end{eqnarray}
	It is easy to see that
	\begin{eqnarray}
		&& |z|^2+|w|^2=|\alpha|^2+|\beta|^2, \qquad
		z-w\xi_{2k}=\alpha-\beta, \qquad  z+w\xi_{2k}=-\overline{\alpha}-\overline{\beta}, \label{Equation:Sec2-zwz-1} \\
		&& |x|^2+|y|^2=|\alpha|^2+|\beta|^2, \qquad
		x-y=\alpha-\beta, \qquad  x+y=-\overline{\alpha}-\overline{\beta}. \label{Equation:Sec2-xyx-1}
	\end{eqnarray}
	
	The following theorem extends the known Hadamard matrix constructions
	from binary and quaternary Legendre pairs to the $k$-ary setting.
	
	\begin{theorem}\label{Theorem:2-1} 
		Let $k\geq 2$ be an integer and
		let $(A, B)$ be a $k$-ary Legendre pair of
		length $\ell$.
		
		\textup{I}. Suppose that $k$ is odd and $z, w\in\{1, \xi_{2k}, \xi_{2k}^2, \ldots, \xi_{2k}^{2k-1}\}$, where $z ,w$ are defined by~$(\ref{Definition:zw})$. 
		Then
		$$
		H(2k, 2\ell+2)=
		\left[\begin{tabular}[]{c|c}
			$\begin{array}{cc}
				z             &  w \\
				\overline{w}  &  -\overline{z}
			\end{array}$  & $\begin{array}{cccccc}
				1 & \ldots & 1 &  \quad  1  & \ldots  &  1  \\
				\xi_{2k} & \ldots & \xi_{2k} &  \quad -\xi_{2k}  & \ldots  &  -\xi_{2k}
			\end{array}$ \\
			\hline
			$\begin{array}{cc}
				1 &  \overline{\xi_{2k}} \\
				\vdots & \vdots  \\
				1 &  \overline{\xi_{2k}}
			\end{array}$ & $\begin{array}{cc}
				C(A) \qquad & \quad \qquad C(B)
			\end{array}$  \\
			\hline
			$\begin{array}{cc}
				1 &  -\overline{\xi_{2k}} \\
				\vdots & \vdots  \\
				1 &  -\overline{\xi_{2k}}
			\end{array}$ & $\begin{array}{cc}
				C(\overline{B})^{T} \qquad & \quad \qquad -C(\overline{A})^{T}
			\end{array}$
		\end{tabular}
		\right]
		$$
		is a Butson Hadamard matrix of order $2\ell+2$.
		
		\textup{II}. Suppose that $k$ is even and $x, y\in\{1, \xi_k, \xi_k^2, \ldots, \xi_k^{k-1}\}$, where $x, y$ are defined
		by $(\ref{Definition:xy})$.
		Then
		$$
		H(k, 2\ell+2)=
		\left[\begin{tabular}[]{c|c}
			$\begin{array}{cc}
				x             &  y \\
				\overline{y}  &  -\overline{x}
			\end{array}$  & $\begin{array}{cccccc}
				1 & \ldots & 1 &  \quad  1  & \ldots  &  1  \\
				1 & \ldots & 1 &  \quad -1  & \ldots  &  -1
			\end{array}$ \\
			\hline
			$\begin{array}{cc}
				1 &  1 \\
				\vdots & \vdots  \\
				1 &  1
			\end{array}$ & $\begin{array}{cc}
				C(A) \qquad & \quad \qquad C(B)
			\end{array}$  \\
			\hline
			$\begin{array}{cc}
				1 &  -1 \\
				\vdots & \vdots  \\
				1 &  -1
			\end{array}$ & $\begin{array}{cc}
				C(\overline{B})^{T} \qquad & \quad \qquad -C(\overline{A})^{T}
			\end{array}$
		\end{tabular}
		\right]
		$$
		is a Butson Hadamard matrix of order $2\ell+2$.
	\end{theorem}
	
	\begin{proof}  
		We prove Part I. Put $\zeta=\xi_{2k}$, $P=C(A)$, $Q=C(B)$,
		and let $\mathbf{1}$ denote the all-one row vector of length $\ell$.
		By assumption, all entries of the displayed matrix are $2k$-th
		roots of unity.
		
		By (\ref{Equation:alpha2+beta2=2}) and (\ref{Equation:Sec2-zwz-1}), each of the first two rows has squared norm
		\[
		|z|^2+|w|^2+2\ell=2\ell+2,
		\]
		and their inner product is
		\[
		zw-wz+\ell\overline{\zeta}-\ell\overline{\zeta}=0.
		\]
		Since every row of $P$ and $Q$ has sum $\alpha$ and $\beta$,
		respectively, the inner products of the first row with rows in the
		two lower block rows are
		\[
		z+w\zeta+\overline{\alpha}+\overline{\beta}=0,
		\qquad
		z-w\zeta+\beta-\alpha=0
		\]
		by (\ref{Equation:Sec2-zwz-1}). For the second row the corresponding inner products are
		\[
		\overline{w}-\zeta\overline{z}
		+\zeta(\overline{\alpha}-\overline{\beta})=0,
		\qquad
		\overline{w}+\zeta\overline{z}
		+\zeta(\alpha+\beta)=0,
		\]
		which follow by taking complex conjugates in (\ref{Equation:Sec2-zwz-1}).
		
		For two different rows in the first lower block row, the first two entries
		contribute $2$, while (\ref{Equation:Sec2-CACA+CBCB})   gives contribution $2\ell$ on the
		diagonal and $-2$ off the diagonal from the remaining entries.
		Thus these rows have squared norm $2\ell+2$ and are mutually
		orthogonal. The same holds for the second lower block row, since
		circulant matrices are normal and hence
		\[
		P^*P=PP^*,\qquad Q^*Q=QQ^*.
		\]
		Finally, for rows belonging to different lower block rows, the
		contribution from the first two entries is zero, whereas that from
		the remaining blocks is
		\[
		PQ-QP=0
		\]
		by (\ref{Equation:Sec2-CACB}). Hence all rows are mutually orthogonal and have squared norm
		$2\ell+2$, so the displayed matrix is an $H(2k,2\ell+2)$.
		
		Part II follows in exactly the same way, replacing
		$\zeta,z,w$ by $1,x,y$, respectively, and using (10) in place
		of (\ref{Equation:Sec2-zwz-1}).  
	\end{proof}
	
	\begin{corollary}\label{cor:prime-butson}
		Let $p$ be an odd prime and let $(A,B)$ be a $p$-ary
		Legendre pair of length $\ell\ge2$. Then there exists a Butson
		Hadamard matrix $H(2p,2\ell+2)$.
	\end{corollary}
	
	\begin{proof}
		By Theorem~\ref{thm:prime-balance}, after multiplying $A$ and
		$B$ independently by suitable $p$th roots of unity, we may assume
		$\alpha=\beta=-1$. Then the quantities defined in (\ref{Definition:zw}) become
		$$
		z=\frac12
		(\alpha-\beta-\overline{\alpha}-\overline{\beta})=1
		$$
		and
		$$
		w=\frac{1}{2\xi_{2p}}
		(\beta-\alpha-\overline{\alpha}-\overline{\beta})
		=\xi_{2p}^{-1}.
		$$
		Thus both $z$ and $w$ are $2p$th roots of unity, and
		Theorem \ref{Theorem:2-1} (I) applies.
	\end{proof}
		
For $k=2$ we may assume $\alpha=\beta=1$, so that
$x=y=-1$. For $k=4$ and even $\ell$ we may assume
$\alpha=0$ and $\beta=1+i$, so that $x=-1$ and $y=i$.
For every odd prime $p$, Part I applies automatically by
Corollary~\ref{cor:prime-butson}.

For odd prime $p$, the construction therefore always yields an
$H(2p,2\ell+2)$. 
As mentioned in the introduction, in many cases there exists a ternary Legendre pair of length $\ell\equiv -1\bmod 3$ but no $H(p,2\ell+2)$.

For general composite $k$, it is not clear whether the conditions
on $z,w$ or $x,y$ in Theorem~\ref{Theorem:2-1} are automatically
satisfied.
	

\section{Ternary Legendre pairs}
\label{ternary}

\quad\ \ We now specialize to $k=3$. By Theorem~\ref{Theorem:length}, every ternary
Legendre pair has length
$$
\ell\equiv -1\pmod 3.
$$
By Theorem~\ref{thm:prime-balance}, after multiplying the two components
independently by suitable third roots of unity, we may assume
$$
\alpha=\sum_{j=0}^{\ell-1}a_j=-1, \qquad  \beta=\sum_{j=0}^{\ell-1}b_j=-1.
$$
Writing $\ell=3h-1$, Theorem~\ref{thm:prime-balance} shows that
each normalized component has composition $(h-1,h,h)$ with
respect to the symbols $1,\xi_3,\xi_3^2$.
The construction of Section~\ref{Subsection:m-sequences} gives ternary
Legendre pairs of lengths $\ell=3^r-1$, $r\geq1$, that is,
$2,8,26,80,\ldots$.

For compactness, in the explicit examples below we sometimes write a ternary
sequence as an exponent sequence $[r_0,\ldots,r_{\ell-1}]$, where
$r_j\in\{0,1,2\}$ represents the entry $\xi_3^{r_j}$.

\subsection{Searches using necklaces with fixed content}

\quad\ \ Our searches produced ternary Legendre pairs of lengths
$\ell=5,11,17$ and $23$. We record existence results, without claiming a
complete enumeration of pairs or their equivalence classes. Explicit witnesses
for lengths $5$, $11$, and $17$ are given in Appendix~\ref{app:small-certificates}.
We give details for $\ell=23$.
After normalization, each component has composition $(7,8,8)$. After a cyclic shift, we may fix the first entry to be $1$. Then there are
$$
    \binom{22}{6}\binom{16}{8}
    =74613\cdot12870\approx9.6\times10^8
$$
possible sequences before equivalences are taken into account.
An example of a ternary Legendre pair of length $23$ is
\begin{eqnarray*}
A&=&[0,1,1,2,2,1,0,0,1,2,0,2,1,1,2,0,2,1,2,1,0,0,2],\\
B&=&[0,1,1,1,2,1,2,2,1,1,2,1,1,0,0,0,2,2,0,2,0,2,0].
\end{eqnarray*}
A direct exact computation gives
\[
    PAF(A,s)+PAF(B,s)=-2,
    \qquad 1\leq s\leq 22.
\]

The search can be reduced substantially by generating only necklaces with
fixed content. A \emph{necklace} is the lexicographically smallest member of
a rotation class; see \cite{Sawada:2001}. Algorithms for necklaces with fixed
content are described in \cite{SawadaETAL:2013}. For $\ell=23$ and content
$(7,8,8)$ there are only $137,181,330$ such necklaces. Thus the balance result
of Section~\ref{non-ex} is not only structural: it fixes the content used by
the necklace generator before any correlation computation is performed.

\subsection{Integral PAF parameters and reflection symmetry}\label{subsec:integral-paf}

\quad\ \ For searches restricted to integral periodic autocorrelations, one can impose
additional exact arithmetic constraints.

\begin{proposition}\label{prop:ternary-integral-paf}
Let $(A,B)$ be a normalized ternary Legendre pair of length
$\ell=3h-1$. Suppose that the nonzero PAFs of one, equivalently both, components are integers. Then there are
integers $x_u$, $1\leq u\leq\ell-1$, such that
$$
    PAF(A,u)=-1+3x_u,  \qquad     PAF(B,u)=-1-3x_u.
$$
Moreover,
$$
    |x_u|\leq \left\lfloor\frac{h-1}{2}\right\rfloor, \qquad
    x_{\ell-u}=x_u,  \qquad   \sum_{u=1}^{\ell-1}x_u=0.
$$
More precisely, if
$$
N_r^A(u)=\#\{j:a_j\overline{a_{j+u}}=\xi_3^r\},\qquad
N_r^B(u)=\#\{j:b_j\overline{b_{j+u}}=\xi_3^r\},
\qquad 0\le r\le2,
$$
then
$$
\begin{aligned}
N_0^A(u)&=h-1+2x_u,&\qquad N_1^A(u)=N_2^A(u)&=h-x_u,\\
N_0^B(u)&=h-1-2x_u,& N_1^B(u)=N_2^B(u)&=h+x_u.
\end{aligned}
$$
\end{proposition}

\begin{proof}
Fix a nonzero lag $u$ and suppose $PAF(A,u)=m\in\mathbb Z$. Since
$$
    PAF(A,u)=N_0^A(u)+N_1^A(u)\xi_3+N_2^A(u)\xi_3^2,
$$
the relation $1+\xi_3+\xi_3^2=0$ implies
$N_1^A(u)=N_2^A(u)$. Hence $m\equiv\ell\equiv-1\pmod3$, and therefore
$m=-1+3x_u$ for a unique $x_u\in\mathbb Z$. The Legendre pair identity gives
$PAF(B,u)=-1-3x_u$. Solving the equations for the three difference counts
gives the displayed formulas, and their nonnegativity yields the pointwise
bound.

For every complex sequence, we have $PAF(A,\ell-u)=\overline{PAF(A,u)}$.
Since the PAF values here are integers, they are real, and hence
$x_{\ell-u}=x_u$. Finally,
$$
    \sum_{u=0}^{\ell-1}PAF(A,u)=|\alpha|^2=1.
$$
Using $PAF(A,0)=\ell$ and $PAF(A,u)=-1+3x_u$ for $u\ne0$ gives
$\sum\limits_{u=1}^{\ell-1}x_u=0$.
\end{proof}

If $\ell$ is odd, the reflection symmetry immediately gives
$\sum\limits_{u=1}^{(\ell-1)/2}x_u=0$.
Thus an integral-PAF search at odd length needs only half of the PAF
coordinates, together with a zero-sum condition.

As an illustration, a ternary Legendre pair of length $29$ is
\begin{align*}
A={}&[0,0,1,0,1,2,1,2,0,2,2,1,1,2,0,\\
   &\phantom{[}2,1,1,2,2,0,2,1,2,1,0,1,0,0],\\
B={}&[1,2,1,1,2,0,0,0,2,1,1,0,2,2,0,\\
   &\phantom{[}2,2,0,1,1,2,0,0,0,2,1,1,2,1].
\end{align*}
For this pair, exact computation gives
\[
\begin{aligned}
(x_1,\ldots,x_{28})={}&(-1,2,-1,0,2,0,-1,0,0,0,0,-1,0,0,\\
&0,0,-1,0,0,0,0,-1,0,2,0,-1,2,-1).
\end{aligned}
\]
In particular, $x_{29-u}=x_u$ and
$\sum\limits_{u=1}^{14}x_u=0$. This provides a convenient exact certificate for the
integral PAF profile and illustrates the reduction available in this subclass.

\subsection{Half-frequency restrictions and compressed searches}
\label{subsec:half-frequency}

\quad\ \ We next consider even lengths. The discrete Fourier transform and the power
spectral density of a sequence $C=[c_0,\ldots,c_{\ell-1}]$ are
$$
    DFT(C,s)=\sum_{j=0}^{\ell-1}c_j\xi_\ell^{js}, \qquad
    PSD(C,s)=|DFT(C,s)|^2.
$$
For a Legendre pair, we have $PSD(A,s)+PSD(B,s)=2\ell+2$, $1\leq s\leq\ell-1$.
Consequently, the usual PSD test discards any candidate component $C$ for
which $PSD(C,s)>2\ell+2$ at some nonzero frequency.

For a divisor $t$ of $\ell$, the compression of
$C=[c_0,\ldots,c_{\ell-1}]$ by the factor $\ell/t$ is the length-$t$ sequence
$C_t=[c'_0,\ldots,c'_{t-1}]$, where
$$
    c'_j=\sum_{s=0}^{\ell/t-1}c_{st+j}.
$$
When $t=2$, write
$$
    C_2=[\gamma_0,\gamma_1], \qquad
    \gamma_0=\sum_jc_{2j}, \quad
    \gamma_1=\sum_jc_{2j+1}.
$$
Then
$$
    PSD(C_2,1)=PSD\left(C,\frac\ell2\right)
    =|\gamma_0-\gamma_1|^2.
$$
Thus compression to length $2$ is exactly the half-frequency test.

The fixed composition from Section~\ref{non-ex} yields a sharper form of this
test. Since an admissible even ternary length has the form
$\ell=6m+2$, each normalized component has composition
$(2m,2m+1,2m+1)$.
For a component $C$, let $e_r$ and $o_r$ be the numbers of occurrences of
$\xi_3^r$ in the even and odd positions, respectively, and put
$$
    d_r=e_r-o_r,\qquad r=0,1,2.
$$

\begin{theorem}\label{thm:half-frequency}
Let $(A,B)$ be a normalized ternary Legendre pair of even length
$\ell=6m+2$. For each component $C\in\{A,B\}$, we have
$$
    PSD\left(C,\frac\ell2\right)=3Q_C,
    \qquad  Q_C=d_0^2+d_0d_1+d_1^2.
$$
The integer $Q_C$ is positive and odd, and $Q_A+Q_B=4m+2$.
\end{theorem}

\begin{proof}
There are equally many even and odd positions, so
$d_0+d_1+d_2=0$. The fixed composition gives
$$
    d_0\equiv0\pmod2,  \qquad d_1\equiv d_2\equiv1\pmod2.
$$
Moreover,
$$
    DFT\left(C,\frac\ell2\right)=d_0+d_1\xi_3+d_2\xi_3^2.
$$
Using $d_2=-d_0-d_1$ and $1+\xi_3+\xi_3^2=0$, a direct calculation gives
$$
    \left|d_0+d_1\xi_3+d_2\xi_3^2\right|^2=3(d_0^2+d_0d_1+d_1^2).
$$
Since $d_1$ is odd, $Q_C>0$, and the parity conditions show that $Q_C$ is
odd. Finally,
$$
    PSD\left(A,\frac\ell2\right)+PSD\left(B,\frac\ell2\right)=2\ell+2=12m+6,
$$
which yields $Q_A+Q_B=4m+2$.
\end{proof}

\begin{corollary}
Under the assumptions of Theorem~\ref{thm:half-frequency}, we have
$$
    PSD\left(A,\frac\ell2\right)
    \equiv
    PSD\left(B,\frac\ell2\right)
    \equiv3\pmod6.
$$
Moreover, in the prime factorization of each half-frequency PSD,
every prime congruent to $2$ modulo $3$ occurs with even exponent.
\end{corollary}

\begin{proof}
The congruence follows because $Q_A$ and $Q_B$ are odd. Also
$$
    Q_C=d_0^2+d_0d_1+d_1^2
       =N_{\mathbb Q(\xi_3)/\mathbb Q}(d_0-d_1\xi_3).
$$
Every prime congruent to $2$ modulo $3$ occurs with even exponent in an
Eisenstein norm; see, for example, \cite[Theorem 3.1]{Brock1988} or \cite{t65}. Multiplication by
$3$ does not affect this assertion.
\end{proof}

The theorem fits directly into the existing compression-decompression search.
Indeed, once an admissible pair $(Q_A,Q_B)$ has been chosen, a vector
$(d_0,d_1,d_2)$ determines the parity content of a component uniquely:
$$
\begin{array}{lll}
 e_0=m+d_0/2, &\quad& o_0=m-d_0/2,\\[1mm]
 e_1=(2m+1+d_1)/2, && o_1=(2m+1-d_1)/2,\\[1mm]
 e_2=(2m+1+d_2)/2, && o_2=(2m+1-d_2)/2.
\end{array}
$$
Equivalently, in terms of the even-position counts,
$$
    PSD\left(C,\frac\ell2\right)=9\bigl(2m+1-(e_1+e_2)\bigr)^2+3(e_1-e_2)^2.
$$
Thus the exact half-frequency restriction is applied first, the corresponding
parity contents are generated next, and only then are decompression and the
full PSD test performed. The original search method is therefore retained,
but its candidate space is reduced before the expensive stages.

This approach yields ternary Legendre pairs of lengths $14$ and $20$.
An explicit length-$14$ witness is included in
Appendix~\ref{app:small-certificates}. For
$\ell=20$ we have $m=3$ and hence $Q_A+Q_B=14$. The admissible complementary
Eisenstein norms give, up to interchanging the components,
$$
    \left(PSD(A,10),PSD(B,10)\right)\in\{(21,21),(3,39)\}.
$$
For example, in the case $(21,21)$ one may take parity contents
$$
    (e_0,e_1,e_2;o_0,o_1,o_2)=(4,4,2;2,3,5)
$$
for each component. This leaves
$$
    \binom{10}{4}\binom{6}{4}
    \binom{10}{2}\binom{8}{3}
    =7.938\times10^6
$$
candidate sequences for each component before the remaining filters.
A search within this candidate class yields the following ternary Legendre pair:
\begin{eqnarray*}
A&=&[0,1,0,0,2,2,0,1,1,2,1,2,1,2,1,1,0,2,2,0],\\
B&=&[0,1,0,2,1,1,1,2,0,1,0,0,2,2,1,2,1,2,2,0].
\end{eqnarray*}
An exact computation verifies $PAF(A,s)+PAF(B,s)=-2$ for every nonzero lag.

\subsection{The length-$32$ search}
\label{subsec:length32}

\quad\ \ The exact half-frequency restriction is particularly effective at length
$32$. Here $m=5$, so $Q_A+Q_B=22$. Since
$$
    Q_C=\left(d_1+\frac{d_0}{2}\right)^2+\frac34d_0^2,
$$
and $d_0$ is even, the bound $Q_C\leq21$ gives
$d_0\in\{0,\pm2,\pm4\}$. A short enumeration yields
$$
    Q_C\in\{1,3,7,9,13,19,21\}.
$$
Consequently, up to interchanging $A$ and $B$, the only possible
half-frequency pairs are
$$
   \left(PSD(A,16),PSD(B,16)\right)\in\{(3,63),(9,57),(27,39)\}.
$$
Thus every ternary Legendre pair of length $32$ must belong to one of these
three spectral cases, which can therefore be searched separately. These are
necessary conditions; no realization of every spectral case is asserted.

A targeted search in the case $(9,57)$ yields the following ternary Legendre pair:
\begin{eqnarray*}
A&=&[0,1,0,0,0,0,1,0,2,2,1,1,2,1,2,1,
      0,1,1,0,2,2,1,2,0,1,0,2,1,2,2,2],\\
B&=&[0,0,1,2,1,0,0,0,2,2,2,2,0,0,1,2,
      0,1,1,2,1,2,1,1,0,2,1,1,2,2,0,1].
\end{eqnarray*}
The even/odd symbol counts are
$$
    A:(6,5,5;4,6,6), \qquad   B:(6,7,3;4,4,8),
$$
so the half-frequency PSDs are indeed $9$ and $57$, respectively.
For $s\ne0$, define
$$
N_r(s)=\#\{j:a_j\overline{a_{j+s}}=\xi_3^r\}+\#\{j:b_j\overline{b_{j+s}}=\xi_3^r\}.
$$
Exact computation gives
$$
    (N_0(s),N_1(s),N_2(s))=(20,22,22), \qquad 1\leq s\leq31.
$$
Hence
$$
    PAF(A,s)+PAF(B,s)=20+22\xi_3+22\xi_3^2=-2
$$
for every nonzero lag, proving that the displayed sequences form a ternary
Legendre pair of length $32$.

The same half-frequency calculation can be used prospectively at larger even
lengths. For example, at $\ell=38$ it gives, up to interchanging the
components,
$$
    \left(PSD(A,19),PSD(B,19)\right)
    \in\{(3,75),(21,57),(39,39)\},
$$
which provides three exact cases for a targeted search.
Here the admissible values with $Q_C\le25$ are
$$
1,3,7,9,13,19,21,25,
$$
and the complementary pairs summing to $26$ are
$(1,25)$, $(7,19)$ and $(13,13)$.
\begin{conjecture}
Ternary Legendre pairs exist for every admissible length $\ell\equiv -1\pmod 3$.
\end{conjecture}
The preceding examples illustrate the role of the arithmetic
restrictions in the search. The balance theorem fixes the symbol
composition, the half-frequency result restricts the admissible
spectral cases at even lengths, and the integral-PAF symmetry
provides an additional reduction when the autocorrelations are
integral. Together with fixed-content necklaces, compression and
the PSD test, these conditions provide a layered search strategy
for larger ternary lengths. The next open cases are $35$ and $38$. 
A broader arithmetic and cyclotomic framework for the structure
of Legendre pairs, motivated in part by the computational
phenomena observed here, will be developed elsewhere.

\appendix
\section{Explicit small-length certificates}
\label{app:small-certificates}

For completeness, we give directly checkable witnesses for the small lengths
stated in Section~\ref{ternary} but not displayed there. In the table below,
an exponent word $r_0\cdots r_{\ell-1}$ denotes the sequence
$[\xi_3^{r_0},\ldots,\xi_3^{r_{\ell-1}}]$. Each row gives one existence
certificate; no total number of solutions or equivalence classes is asserted.
The length-$11$ and length-$14$ witnesses were independently reconstructed
when preparing these certificates and are not a reproduction of the original
search logs.

\begin{center}
\renewcommand{\arraystretch}{1.2}
\begin{tabular}{c@{\qquad}l@{\qquad}l}
\hline
$\ell$ & Exponent word for $A$ & Exponent word for $B$ \\
\hline
5 & \texttt{12021} & \texttt{12021} \\
11 & \texttt{00120211212} & \texttt{00021122121} \\
14 & \texttt{00120202121121} & \texttt{00112211022201} \\
17 & \texttt{01212201010221210} & \texttt{12110220002201121} \\
\hline
\end{tabular}
\end{center}

For each listed length, write $\ell=3h-1$. Both components have composition
$(h-1,h,h)$. For every nonzero lag $u$, define the combined quotient counts
\[
M_r(u)=\#\{j:a_j\overline{a_{j+u}}=\xi_3^r\}
       +\#\{j:b_j\overline{b_{j+u}}=\xi_3^r\},
\qquad 0\le r\le2.
\]
Direct integer evaluation at every $1\le u\le\ell-1$ gives
\[
\begin{array}{c|c}
\ell &(M_0(u),M_1(u),M_2(u))\text{ for every nonzero }u\\
\hline
5 &(2,4,4)\\
11&(6,8,8)\\
14&(8,10,10)\\
17&(10,12,12).
\end{array}
\]
In each case this is $(2h-2,2h,2h)$, so
\[
PAF(A,u)+PAF(B,u)
=(2h-2)+2h\xi_3+2h\xi_3^2=-2.
\]
This verifies the Legendre pair condition without floating-point arithmetic.
The accompanying file \texttt{verify\_existence.py} performs these checks,
checks all explicit pairs in the main text, and generates normalized
$m$-sequence certificates for lengths $2$, $8$, and $26$ using the construction
of Section~\ref{Subsection:m-sequences}. Thus every nontrivial admissible
length $2\le\ell\le32$ is covered either by that construction or by an
explicit checked pair. The verification establishes existence, not the
completeness or original output count of any search.

\section*{Acknowledgments}
	
	The authors thank the organizers and referees of SETA 2026 for helpful comments and discussions. This work was developed in part during a visit of Huaning Liu to RICAM, Linz, and during a visit to the Hangzhou International Innovation Institute of Beihang University (H3I), Sino-French Laboratory for Mathematics (SFLM), Hangzhou, P. R. China.
		The second author was supported by National Natural Science Foundation of China under Grant No. 12571008,
		and the Science and Technology Program of
		Shaanxi Province of China under Grant No. 2024JC-JCQN-04.
		The third author was funded in part by the Austrian Science Fund (FWF) (Grant No.\ 10.55776/PAT4719224).

\end{document}